%% file: main.tex
\documentclass[letterpaper, 10 pt, conference]{cssconf}

\IEEEoverridecommandlockouts
\usepackage{lipsum}
\usepackage{xcolor}
\usepackage{cite}
\usepackage{graphicx}
\usepackage{float}
\usepackage{tikz}
\usepackage{circuitikz}
\usetikzlibrary{shapes.geometric, arrows, decorations.markings, arrows.meta}
\usepackage{amsmath}
\usepackage{amssymb,amsfonts,amscd}
\input{math.tex}

\usepackage{layouts}

\title{\LARGE \bf
Kalman-based approaches for online estimation of bioreactor dynamics from
fluorescent reporter measurements$^*$
}

\author{Rand Asswad$^{1,\dagger}$, Eugenio Cinquemani$^{1}$, Jean-Luc Gouzé$^{2}$%
\thanks{$^*$Work supported in part by project Ctrl-AB [ANR-20-CE45-0014]}
\thanks{$^{1}$Université Grenoble Alpes, Inria, Grenoble, France}
\thanks{$^{2}$Université Côte d'Azur, Inria, INRAE, CNRS, Sorbonne Université,
        Macbes Team, Sophia Antipolis, France}
\thanks{$^\dagger$Corresponding author, {\tt\small rand.asswad@inria.fr}}
}

\begin{document}

\maketitle
\thispagestyle{empty}
\pagestyle{empty}


\begin{abstract}

We address online estimation of microbial growth dynamics in bioreactors
from measurements of a fluorescent reporter protein
synthesized along with microbial growth.
We consider an extended version of standard growth models
that accounts for the dynamics of reporter synthesis.
We develop state estimation from sampled, noisy measurements
in the cases of known and unknown growth rate functions.
Leveraging conservation laws and regularized estimation techniques,
we reduce these nonlinear estimation problems to linear
time-varying ones, and solve them via Kalman filtering.
We establish convergence results in absence of noise and show performance
on noisy data in simulation.

\end{abstract}

\section{Introduction}\label{intro}

Estimation and control of cellular growth in bioreactors have been dedicated a significant amount of work over several decades \cite{bastin_1990}.
Nowadays, these challenges receive renewed attention
in relation with new biotechnological
developments. Among the new research frontiers
is the feedback control of microbial consortia \cite{clementschitsch_2006}. Motivated by the potential to outperform single species in biosynthetic
processes and other applications \cite{bensaid_2017}, control of microbial consortia requires in the first place to be able to discern single species dynamics by suitable real-time monitoring, \textit{e.g.} the use of fluorescent
reporter proteins expressed along with the growth of the different species.

In a recent paper \cite{alex_2021}, we have explored 
(local) observability of the prototypical microbial consortium model from \cite{mauri_2020} under different assumptions about the observed variables. The construction and analysis of observers for (noisy) time-sampled measurements (as encountered frequently, \textit{e.g.} with the use of pipetting robots \cite{bertaux_2022})
have not been explored in detail. Given the nonlinear nature of these models, the problem is
nontrivial. General observer synthesis approaches exist 
\cite{besancon_2007,gauthier_2001}, however their viability depends on perfect knowledge of
the system, which hardly applies
to biological systems, and/or on idealized assumptions about the measurement process.
Robust estimation approaches can tackle in part modelling uncertainty \cite{hogg_1979}, but they are generally not conceived for sampled data. Nonlinear generalizations of Kalman filtering exist that cope explicitly
with sampled noisy data \cite{jazwinski_1970,sarkka_2007}, but performance is usually not guaranteed even with a perfectly known model.

In this paper, we focus on the construction of state estimators for a single species growing in a bioreactor from indirect, sampled
fluorescent reporter measurements of the species abundance. We consider classical nonlinear growth
models \cite{bastin_1990,smith_1995})
modified to account for the dynamical relationship between 
fluorescence abundance
and biomass synthesis. 
As explained above, besides its interest per se, the problem constitutes a key step toward estimation of dynamics of
microbial communities. Much literature has been dedicated to the problem in the case of
direct observations of biomass abundance
\cite{bastin_1990,neeleman_2001}
or gaseous outflow \cite{didi_2023}.
None of these approaches applies to the augmented model we consider with sampled, noisy measurements.

We consider two cases, perfect knowledge of the functional form of the growth law, and lack of this knowledge. 
We apply hybrid (continuous-time dynamics, discrete-time measurements) Kalman filtering to linear, time-varying versions of the original nonlinear estimation problems, obtained by treating nonlinearities as unknown time-varying inputs. In the first case, we show that the unknown input can be pre-estimated thanks to the system conservation laws. In the second case, we use a Bayesian regularization approach to linearly estimate the unknown input
along with the other system state variables. We demonstrate in simulation that our approaches 
perform well in presence of noise, 
and we prove that the estimator
for a known growth law is an observer, \emph{i.e.}, it enjoys deterministic convergence in absence of noise.

The paper is organized as follows. In Sec. \ref{model}, we describe bioreactor growth models and illustrate the properties of our model. In Sec. \ref{observer} we review concepts of observability
and Kalman filtering estimation that are used in the sequel.
In Sec. \ref{known} we introduce our estimation method 
for a perfectly known growth law,
and prove its convergence properties.
In Sec. \ref{unknown} we present our method to
address state estimation for unknown growth laws.
In Sec. \ref{results}, based on noisy simulated data, we compare performance
of our methods with the reference approach of Extended Kalman filtering.
Conclusions and perspectives are drawn in Sec. \ref{conclusion}.

\section{Model and properties}\label{model}

The model of a continuously stirred-tank bioreactor (CSTR)
for microbial growth in constant volume
without accounting for gaseous outflow is generally expressed as
\begin{equation}
\dxt = N r(x,t) + d(\Xin - x) 
\label{eq:chemostat}
\end{equation}
where
$x(t)\in\Omega=\R_+^n$ represents the reactants concentrations vector,
$\Xin\in\R_+^n$ is the reactor feed concentrations vector,
$N\in\R^{n\times q}$ is the stoichiometric matrix,
$r(x,t)\in\R_+^q$ is the reaction rates vector ($\gL\h$),
and $d\in\R_+^*$ is the dilution rate ($\h$) \cite{bastin_1990}.
Concentrations (in $\gL$) are relative to the total culture
volume.
We assume a strictly positive dilution rate $d>0$.
Consider in particular the dynamical system describing
the growth of \textit{Escherichia coli} (\textit{E. coli}) bacteria 
along with synthesis of a fluorescent protein. 
We assume as in \cite{alex_2021} that the latter
takes part of the resources away from the synthesis of the self-replicative
biomass. For $e(t)$, $f(t)$ and $s(t)$ the biomass, fluorescent protein, and substrate concentration, in the same order, and $x(t)=\transp{(s(t),e(t),f(t))}\in\R_+^3$, the model becomes~\cite{hidde_2010,alex_2021,pavlou_2022}
\begin{equation}
\begin{aligned}
    \dst &= &-\frac{1}{\gamma}&\mu(x,t)e +d(\Sin - s)\\
    \det &= &(1-\alpha)&\mu(x,t)e -de\\
    \dft &= &\alpha&\mu(x,t)e -df
\end{aligned}
\label{eq:fluo}
\end{equation}
where $\mu(x,t)\in\R_+$, the growth rate per capita
at time $t\geq t_0$, 
may depend on time due to quantities that are not part of the model state (pH, temperature) \cite{bastin_1990}.
Constant $0<\alpha<1$ is the proportion of substrate import dedicated to
fluorescent protein synthesis, while $\Sin$ represents the inflow substrate concentration.
Overall, this is
a constant yield model where the substrate
uptake rate per capita is given as
$\mu(x,t)/\gamma$ with $\gamma$ the 
growth yield coefficient.
The values of $\alpha$, $\Sin$, and $\gamma$ are assumed known and strictly
positive.
The model is in the form of (\ref{eq:chemostat})
with $N=\transp{(-1/\gamma,1-\alpha,\alpha)}\in\R^{3}$,
$r(x,t)=\mu(x,t)e(t)\in\R_+$, and $\Xin=\transp{(\Sin,0,0)}\in\R_+^3$.

For the problem of online state observation, in section \ref{known}
we consider the specific growth rate known and strictly a function of the
system's state.
In section \ref{unknown}, we will instead address the case where the functional
form of $\mu(\cdot)$ is unknown, and thus treat it as an unknown function of
time.

\subsection{Conservation relations}

In a biochemical system, linear dependence of the rows of $N$ defines conservation relations.
The rows of $N$ are linearly dependent if and only if
the left null-space of $N$ has
dimension $\ell = n - \rk(N)>0$,
which corresponds to the number of conservation relations.
A conservation matrix $\Pi\in\R^{\ell\times n}$ is defined as a matrix 
with rows forming a basis of the left-null space of $N$.
Consequently, $\Pi N=0$ and $\rk(\Pi) = \ell$.
The conservation matrix, which is not unique
\cite{fjeld_1974,heinrich96_sto},
projects the system's state into an auxiliary
state $\sigma = \Pi x$ of dimension $\ell$.
We also refer to this definition as the conservation law.
From (\ref{eq:chemostat}),
\begin{equation}
\dsigt = d(\Pi\Xin - \sigma)
\label{eq:conservation_dyn}
\end{equation}
(in the ``batch mode'' $d\equiv0$, $\dsigt=0$, whence the name ``conservation'').
The conservation manifold (nullcline of the conservation law) is defined as $\Sigma=\cb{x\in\R^n,\Pi\dxt=0}$.
From \eqref{eq:conservation_dyn}, 
$\sigma(t) = \Pi\Xin + e^{-d(t-t_0)}(\sigma(t_0) - \Pi\Xin)$,
which shows that 
$\Sigma$ is positively invariant.
Since $d>0$, $\Sigma$ is also exponentially attracting, thus it  
is asymptotic invariant \cite{fjeld_1974}.

System \eqref{eq:fluo} has dimension $n=3$  and $\rk(N)=1$. Hence
there are $\ell=2$ conservation relations. Notably, we choose
\begin{equation}\label{eq:cons_fluo}
\Pi = \pmat{
1 & 0 & \dfrac{1}{\alpha\gamma}\\
0 & 1 & -\dfrac{1-\alpha}{\alpha}
} \implies\left\{
\begin{aligned}
    \sigma_1 &= s + \frac{1}{\alpha\gamma}f,\\
    \sigma_2 &= e - \frac{1-\alpha}{\alpha}f.
\end{aligned}\right.
\end{equation}

\subsection{Steady states}\label{equilibria}

Regardless of the form of $\mu(x,t)$, system
(\ref{eq:fluo}) always admits an equilibrium
$x_\text{washout} = \transp{(\Sin,0,0)}$, which is globally asymptotically
stable on $\Omega$ if the set of equilibria is a singleton and unstable
otherwise.
The existence and stability of other equilibria depends on $\mu(x,t)$
and on the system parameters, in this case $d$, $\alpha$, and $\Sin$.

We focus on the  Monod law \cite{bastin_1990}. In this case $\mu$
takes the form
\begin{equation}\label{eq:monod}
    \mu(s) = \mumax \frac{s}{k_s + s}, \quad s\geq 0,
\end{equation}
with $k_s>0$ known as half-saturation constant and
$\mumax>0$ the maximal growth rate.
Eq.~\eqref{eq:monod} is strictly increasing,
therefore it defines a bijection from $[0,\infty)$ to $[0,\mumax)$.
If ${d/(1-\alpha) < \mu(\Sin)}$, \eqref{eq:fluo} admits another equilibrium
$x^*=\transp{(s^*, \gamma(1-\alpha)(\Sin-s^*), \gamma\alpha(\Sin-s^*))}$
with ${s^*=\mu^{-1}\pp{d/(1-\alpha)}}$.
If this exists, it is globally asymptotically stable
on any open subset of $\Omega\setminus\cb{e=0}$.
From now on we work under this assumption.
These results extend similar known results for the subsystem
$(s,e)$ without fluorescence synthesis \cite{smith_1995}.
They apply similarly to several other common growth models,
such as Teissier law or Contois law \cite{bastin_1990}.

\section{Observability and state estimation} \label{observer}

We review the system's observability and the definition of a state
observer in the idealized case of noiseless measurements in continuous time.
A dynamical system describes the time evolution of the system's state
$x(t)\in\Omega\subset\R^n$. Typically, information about the state variables
is obtained through measurements $y(t)=h(x(t))\in\R^p$
(often referred to as the system's output),
where $h$ is a function assumed $\CC^\infty(\Omega)$.
The problem of interest in this paper is the construction of a state
observer $\obsx(t)\in\R^n$ from the measurements $y(t)$ for all $t\geq t_0$.
$\obsx(t)$ is an observer if
(i) $\norm{x(t)-\obsx(t)}$ converges to zero as $t$ tends to $\infty$,
and (ii) $\obsx(t)=x(t)\implies\obsx(t+\tau)=x(t+\tau),\forall\tau\geq0$.
If the first condition holds for any $x(t_0)$ and $\obsx(t_0)$, the observer
is \textit{global}.
If the convergence rate for (i) can be tuned, the observer is \textit{tunable}.
A tunable observer can be constructed if the system is \textit{observable},
which is discussed in subsection \ref{observability} \cite{besancon_2007}.
An open loop (non-tunable) observer (\textit{i.e.} with no output feedback)
is also referred to as a \textit{detector}.

The available measurements for (\ref{eq:fluo}) are
fluorescent protein concentrations $y(t)=f(t)$, thus
$h$ is defined by 
${y(t)=C x(t)}$ with $C=(0,0,1)$. 

\subsection{Observability}\label{observability}

The \textit{local observability} of a dynamical system $\dxt=f(x)$ is studied
via its observation space $\OO(h)$, that is the vector field containing
$h$ and its Lie derivatives along the vector field $f$.
The Lie derivative along $f$ of a function $\phi:\R^n\to\R^p$
is $L_f(\phi) = \partials{\phi}{x}f$.
The system is locally observable if and only if there exists $n$ elements
of $\OO(h)$ defining a diffeomorphism around $x$ for all $x\in\Omega$
\cite{besancon_2007}.
By the chain rule, $\dyt(t)=\partials{h}{x}(x(t))\dxt(t)=L_f(h)(x(t))$.
Ergo, the observation space is the space of time-derivatives of the outputs.

For (\ref{eq:fluo}), $y=f$, therefore $\dyt = \alpha \mu(x,t)e - d\cdot y$.
Setting $z = (d\cdot y + \dyt)/\alpha = \mu(x_g,t)e$ defines a reduced
dynamical system for the state $x_g=\transp{(s,e)}$ with 
output
${z=h_g(x_g)=\mu(x_g,t)e}$, provided $\mu(x,t)$ does not involve $f$.
Consequently, $\OO(h_g)\subset\OO(h)$, and since a diffeomorphism exists
from $f$ to $z$ by construction, the system is locally observable if the
reduced system is observable.

Biologically, output $z=\mu(x_g,t)e$ corresponds to measurements
of biogas flow rate in the anaerobic digestion process, which is studied and
proven locally observable under suitable assumptions in \cite{didi_2023}.
Observability of (\ref{eq:fluo}) for the Monod growth law follows from
the results in \cite{didi_2023}.

\subsection{Estimation problem and Kalman filtering}\label{kalman}

We stated the observer problem in continuous time, but in truth 
measurements are in discrete time.
This is not an issue in practice if the time steps are sufficiently small
(for example, for continuously operating sensors built on the bioreactor). 
In many experimental scenarios of interest, however, measurements are obtained
by time-consuming procedures (see \textit{e.g.} \cite{bertaux_2022})), hence 
motivates us to explicitly account for the sampled nature of the data.
Moreover, measurements are typically noisy due to biophysical uncertainty in
the measurement process 
This gives rise to an \textit{estimation problem}, as follows. 

Let measurements at times $t_0<t_1<\cdots<t_{M-1}<t_M$ 
be given by
$y_k = C x(t_k) + v_k$,
where $v_k$ is an independent Gaussian noise with zero mean and known
covariance matrix $R_k\in\R^{p\times p}$ (scalar variance if $p=1$).
The problem is to estimate the continuous-time state $x(t)$ from the discrete-time measurements $y_k$. It 
can be addressed by a hybrid (continuous-discrete) Kalman filter \cite{jazwinski_1970}.
This filter is optimal for linear state-space models, while our system dynamics
are nonlinear. Though variants exist that cope with nonlinear
dynamics, such as the Extended KF (EKF) or the Unscented KF (UKF)
\cite{jazwinski_1970,anderson_1979,sarkka_2007},
convergence and performance of these filters are hardly proven.
We instead propose to reformulate the estimation problem by treating 
the nonlinearities of the system dynamics as a time-varying input. This input
is either pre-estimated (as explained in Sec. \ref{known}) or described
by a probabilistic prior in the form of a conveniently chosen linear stochastic
process (as explained in Sec. \ref{unknown}).
In both cases, this leads to a hybrid filtering problem on a linear
(time-varying) stochastic system of the form
\begin{equation}\label{eq:kalman}
\begin{aligned}
    \dxt(t) &= A(t)x(t) + b(t) + G(t) w(t),\\
    y_k &= C_k x(t_k) + v_k,
\end{aligned}
\end{equation}
with 
$A(t)\in\R^{n\times n},b(t)\in\R^n,C_k\in\R^{p\times n},G(t)\in\R^{n\times m},
Q(t)\in\R^{m\times m}$ and $R_k\in\R^{p\times p}$ known. The process noise 
$\{w(t),t\geq t_0\}$ is formally the derivative of the Brownian motion
\cite{jazwinski_1970}, \textit{i.e.} it is a zero-mean, stochastic process with
$\EE{w(t)\transp{w(\tau)}}=Q(t)\delta(t-\tau)$, $\delta$ being the Dirac distribution.

The hybrid KF, which is optimal for this class of models,
alternates ``prediction'' based on the dynamical model and ``update'' based on a new measurement.
The prediction step consists of integrating the 
the system's dynamics for the state estimate $\obsx_{k|k-1}(t)$ and a Riccati
equation for the state covariance matrix $P_{k|k-1}(t)$
over each time interval $[t_{k-1},t_k]$,
with $\obsx_{k|k-1}(t_{k-1})=\obsx_{k-1|k-1}$ and
$P_{k|k-1}(t_{k-1})=P_{k-1|k-1}$.
The update step applied to $\obsx_{k|k-1}(t_k)$ and
$P_{k|k-1}(t_k)$ for a new measurement $y_k$  
is identical to the discrete KF \cite{jazwinski_1970}. 

\subsection{Parameter tuning from data}\label{ML}

In our estimator designs, Eq. (\ref{eq:kalman}) will comprise unknown parameters
$\theta\in\Theta$. We propose (and later demonstrate in simulation) to estimate
these parameters from a preliminary data set $Y_M=\cb{y_0,\ldots,y_M}$
through a Maximum Likelihood (ML) approach,
which closely relates with the filter's prediction error methods
for tuning filters \cite{lutkepohl_1991a}.

Given $Y_M$, ML estimation is the problem of maximizing
the value of the joint distribution
$f(Y_M|\theta)$ over $\theta\in\Theta$.
From the Bayes law,
$f(y_M, Y_{M-1}|\theta) = f(y_M|Y_{M-1},\theta)f(Y_{M-1}|\theta)$.
By recursive application one can write
$f(Y_M|\theta)$ as the product of the conditional
(Gaussian) densities $f(y_k|Y_{k-1},\theta)$ for $k=0,\ldots,M$.
Mean and covariance of these densities correspond to the Kalman prediction
$C_k\obsx_{k|k-1}$ and its covariance $S_k=C_k P_{k|k-1} \transp{C_k} + R_k$
\cite{lutkepohl_1991a,pavlou_2022}.
Writing ML as the equivalent negative log-likelihood minimization and working
out the equations, the ML estimator of $\theta$ is
\small
\begin{equation}\label{eq:ML_pb}
\begin{aligned}
\obstheta &= \argmax\limits_{\theta\in\Theta} f(Y_M|\theta)\\
    &= \argmin\limits_{\theta\in\Theta}
        \bb{-\sum\limits_{k=0}^{M} \log f(y_k|Y_{k-1},\theta)}\\
    &= \argmin\limits_{\theta\in\Theta}
        \bb{-\frac{1}{2}\sum\limits_{k=0}^{M} p\log(2\pi) + \log\abs{S_k}
            + \transp{\erry_{k}}S_k^{-1}\erry_{k}}.
\end{aligned}
\end{equation}
\normalsize
For any candidate $\theta$, all of these quantities can be
computed by the corresponding Kalman filtering iteration, which enables fast
numerical solution of the optimization problem.

\section{Approach for known $\mu(x)$} \label{known}

In this section, leveraging conservation laws, we propose a pipeline for state
estimation for the case where $\mu$ is a known function of the state variables.
We develop this for system (\ref{eq:fluo}) with Monod growth law, yet,
the approach is viable for other growth laws and similar 
biochemical systems.
The resulting estimator is proven to be a global observer.

As seen in section \ref{model}, nonlinearities of CSTR dynamics are
in the reaction rate vector $r(x,t)$ of (\ref{eq:chemostat}).
For system (\ref{eq:fluo}), it is however possible to compute a \textit{pre-estimate}
$\obsmu(t)$ of the specific growth rate $\mu(x,t)$.
Replacing $\mu(x,t)$ in (\ref{eq:fluo}) with $\obsmu(t)$ yields linear
time-varying dynamics,
eligible for linear Kalman filtering.
We next show how to compute $\obsmu(t)$ in the case of the
Monod growth rate $\mu(s)$.

\subsection{Proposed pipeline}

The Monod law is a function of the unknown $s(t)$.
By the conservation law (\ref{eq:cons_fluo}),
in the (exponentially attracting) conservation manifold, $s=\sigma_1 - \frac{1}{\alpha\gamma}f$.
In view of this, we define a pre-estimate $\obss_\sigma$ of $s$ by $\obss_\sigma=\hat \sigma_1 - \frac{1}{\alpha\gamma}\hat y$, with $\hat \sigma_1$ calculated based on the dynamics of $\sigma_1$ and $\hat y$ a suitable smoothening of the measurements $y_k$ of $f$, as detailed next.

The conservation variable $\sigma_1$ is estimated through an open-loop detector
$\obssig_1(t)$ with dynamics as in (\ref{eq:conservation_dyn}), resulting in 
\begin{equation}\label{eq:sig1}
    \obssig_1(t) = \Sin + e^{-d(t-t_0)}(\obssig_1(t_0) - \Sin).
\end{equation}
This converges to $\Sin$ exponentially at a rate $d$.

To calculate $\obsy$, a Kalman-based regularization method is used. 
In the same spirit as \cite{denicolao_1997},  
the method introduces a regularizing linear stochastic dynamical model for the
time profile $f$ to estimate, and solves the problem by linear state-space
estimation \cite{denicolao_2003}.
We consider a second-order linear stochastic differential equation
$\ddot{y}(t) = \beta w(t)$ with $\{w(t)\}$ a scalar white noise process,
and assume that the noisy measurements $y_k$ of $f$ are noisy measurements
of process $y$.
Such model of the data can obviously be expressed in the form of
(\ref{eq:kalman}) as
\begin{equation}\label{eq:presmooth}
    \ddt\bmat{y(t)\\\dyt(t)} =
    \bmat{0 & 1\\0 & 0}\bmat{y(t)\\\dyt(t)} + \bmat{0\\\beta}w(t)
\end{equation}
with $C=(1,0)$.
An optimal parameter $\obs{\beta}$ is the solution to the ML problem discussed
in subsection \ref{ML}. This provides a model of the data-generating process $f$ independent of the knowledge of $f$. Applying a Kalman filter based on this model
to the measurments $y_k$ yields estimates of $y$, \textit{i.e.} of $f$. The chosen form of process $y$ 
is known to penalize overly irregular solutions \cite{denicolao_1997,pavlou_2022}, whereas the ML estimation of $\beta$ ensures an appropriate tradeoff between noisy data interpolation and regularity of the solution. The filter produces 
piecewise-continuous predictions
$\obsy_{k|k-1}(t)$ and filtered estimates $\obsy_{k|k}$. 
Both choices being equally viable as a definition of the estimate $\obsy(t)$ sought,
we simply refer to $\obsy(t)$ from now on and add details only when necessary.

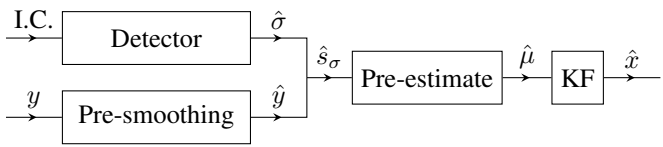
\begin{figure}[H]
\centering
\newcommand{\inArrow}{\tikz\draw[-stealth] (-1pt,0) -- (1pt,0);}
\tikzstyle{block} = [rectangle, minimum height=2em, text centered, draw=black]
\tikzstyle{point} = [shape=coordinate]
\tikzstyle{fixed} = [minimum width=2.5cm]
\begin{tikzpicture}
    \node (s1) at (0,3em)    [point] {};
    \node (s2) at (0,0) [point] {};

    \node (conservation) [block, fixed, right of=s1, xshift=1cm] {Detector};
    \node (presmooth)    [block, fixed, right of=s2, xshift=1cm] {Pre-smoothing};

    \node (m1) at (4cm,3em)      [point] {};
    \node (m)  at (4cm,1.5em) [point] {};
    \node (m2) at (4cm,0)   [point] {};

    \node (preestimate) [block, right of=m, xshift=0.6cm] {Pre-estimate};
    \node (kf) [block, right of=preestimate, xshift=1cm] {KF};
    \node (f) [point, right of=kf, xshift=0.1cm] {};

    \draw (s1)
        -- node[anchor=south]{I.C.} node{\inArrow} (conservation)
        -- node[anchor=south]{$\obssig$} node{\inArrow} (m1)
        -- (m);
    \draw (s2)
        -- node[anchor=south]{$y$} node{\inArrow} (presmooth)
        -- node[anchor=south]{$\obsy$} node{\inArrow} (m2)
        -- (m);
    \draw (m)
        -- node[anchor=south]{$\obss_\sigma$} node{\inArrow} (preestimate)
        -- node[anchor=south]{$\obsmu$} node{\inArrow} (kf)
        -- node[anchor=south]{$\obsx$} node{\inArrow} (f);
\end{tikzpicture}
\vspace{-1em}
\caption{Diagram of the proposed pipeline.}\label{fig:pipeline}
\end{figure}

Armed with $\obssig_1$ and $\obsy$, one can now compute 
$\obsmu(t)=\mu(\obss_\sigma(t))$ 
with $\obss_\sigma(t) = \obssig_1(t) - \frac{1}{\alpha\gamma} \obsy(t)$,
and plug it into \eqref{eq:fluo} to approach Kalman filtering estimation of the whole system state.
Let $\errmu(t) = \mu(s(t)) - \obsmu(t)$ be the pre-estimation error
of $\mu$ at time $t$. Then
\begin{equation}\label{eq:sigma_model}
\begin{aligned}
    \dxt &= N(\obsmu(t) + \errmu(t))e + d(\Xin - x)\\
         &= A_\sigma(t)x + b + N w(t)
\end{aligned}
\end{equation}
with $w(t)=\errmu(t)e(t)$ as the process error term, where 
\begin{equation*}
A_\sigma(t) = \bmat{
    -d & -\obsmu(t)/\gamma      & 0\\
    0  & (1-\alpha)\obsmu(t) -d & 0\\
    0  & \alpha\obsmu(t)        & -d
}, b=\bmat{d{\cdot}\Sin\\0\\0}
\end{equation*}
are now known.
To run the linear Kalman filter of subsection~\ref{kalman} for this system
(with $G=N$ and $C=(0,0,1)$ in Eq. (\ref{eq:kalman})), we exploit the analytic
expression of $\obssig_1(t)$ in (\ref{eq:sig1}), for a given $\obssig_1(t_0)$,
while $\obsy$ and $\obsmu$ can be obtained online sequentially
for every interval $\Tk=[t_k,t_{k+1})$ for $k=0,\ldots,M-1$,
as shown in Figure \ref{fig:pipeline}.

The process noise $w(t)$ remains unknown, and linearly dependent on the state
variable $e$. In the next subsection, an approximate characterization of $w(t)$
is proposed, which is used to tune the KF.
The specific choice of $Q(t)$ is especially important in determining performance
in presence of measurement noise.

\subsection{Tuning of $Q(t)$}

We now derive an approximate expression for $Q(t)$ that characterizes
$w(t)=\errmu(t)e(t)$.
Consider the linear approximation of $\mu$ around $s(t)$.
For sufficiently small $\errs_\sigma=s - \obss_\sigma$, 
\begin{equation}\label{eq:lin_approx}
\mu(\obss_\sigma)=\mu(s - \errs_\sigma)\approx \mu(s) - \mu'(s)\errs_\sigma.
\end{equation}
Then $\errmu(t) = \mu'(s(t))\errs_\sigma(t)$ implies 
$w(t)=\mu'(s(t))e(t)\errs_\sigma(t)$.
We choose to study the process noise around the system equilibrium.
From (\ref{eq:cons_fluo}), (\ref{eq:sig1}) and subsection \ref{equilibria},
it can be shown that for $t\in\Tk, \errs_\sigma(t) = v_k/\alpha\gamma$
at the system's steady states. Consequently,
$w(t) = \frac{\mu'(s^*)e^*}{\alpha\gamma} v_k$ for $t\in\Tk$,
near the system's equilibria.
It follows that $\cb{w(t)}$ is a piecewise constant stochastic process.
By integrating $\EE{w(t)w(\tau)} = Q(t)\delta(t-\tau)$ over $\R$
we obtain the piecewise-constant function
$Q(t)=\pp{\frac{\mu'(s^*)e^*}{\alpha\gamma}}^2 R_k (t_{k+1} - t_k)$ for $t\in \Tk$.

\subsection{Asymptotic stability of the filter}

We now prove that the proposed filtering pipeline defines a \textit{global observer} in the sense of Section \ref{observer}. That is, we consider 
noiseless measurements and look at convergence at the discrete measurement times \cite{besancon_2007}.
We first prove convergence of the observer error, then
its positive invariance, \textit{i.e.}
that $x_k - \obsx_k = 0$ implies $x_{k+l} - \obsx_{k+l} = 0$ for $l>0$.

As seen before, the KF is applied to the pre-estimated system
(\ref{eq:sigma_model}) 
with $A_\sigma(t)x=N\obsmu(t)e-dx$.
In the absence of measurement noise, $\obss_\sigma(t)$ converges to $s^*$
at the same rate as $\obssig_1$ (globally, exponentially at rate $d$).
Consequently, $\obsmu(t)=\mu(\obss_\sigma(t))$ converges to
$\mu^*=\mu(s^*)=d/(1-\alpha)$
and the system (\ref{eq:sigma_model}) is asymptotically time-invariant
with $A_\sigma(t)$ converging globally to
\begin{equation}\nonumber
A = \bmat{
    -d & -\mu^*/\gamma & 0\\
    0 & 0 & 0\\
    0 & \alpha\mu^* & -d
}.
\end{equation}
By integrating the linear time-invariant dynamics over $\Tk$,
\begin{equation}\nonumber
\begin{aligned}
    \dxt(t) &= A x(t) + b + N w(t)\\
    x_{k+1} &= e^{A(t_{k+1}-t_k)} x_k + \int_\Tk e^{A(t_{k+1}-t)} (b + N w(t))\dt.
\end{aligned}
\end{equation}
One can verify that here, $e^{A\tau}N=N,\forall\tau$.
We obtain the difference equation, for all $k$, 
\begin{equation}
x_{k+1} = F_k x_k + b_k + N w_k
\end{equation}
with $F_k=\exp(A(t_{k+1}-t_k))$ the state transition matrix \cite{jazwinski_1970},
$b_k=\pp{\int_\Tk \exp(A(t_{k+1}-t_k))\dt}b$,
and $w_k=\int_\Tk w(t)\dt$.
Let us fix a time step $\Delta t=t_{k+1}-t_k,\forall k$. Let $F=F_k$.
For a fixed $Q_c$ for the continuous system, $\EE{w_k w_l}=Q \delta_{k,l}$
with $Q = Q_c\Delta t$. 

The equivalent discrete Kalman equation is
\begin{equation}\label{eq:kalman_discrete}
\obsx_{k+1|k} = (F - K_k C) \obsx_{k|k-1} + K_k y_k + b_k
\end{equation}
with $K_k\in\R^{n\times 1}$ the Kalman gain.
If the pair $[F,C]$ is completely detectable and $[F,NG]$
is completely stabilizable for any $G$ such that $G\transp{G}=Q$,
then (i) for any positive semidefinite initial matrix $P_0$,
$P_{k+1|k}$ converges to the solution $P_\infty$ of the algebraic
Riccati equation, which defines the associated gain
$K_\infty = F P_\infty\transp{C}(C P_\infty \transp{C} + R)^{-1}$,
(ii) all eigenvalues $\lambda$ of $F-K_\infty C$ verify $\abs{\lambda}<1$
\cite{anderson_1979}. The pair $[F,C]$ is completely detectable if
$\exists K\in\R^3: \abs{\lambda}<1, \forall\lambda\in\Sp(F-KC)$ \cite{anderson_1979}.
Let $K=\transp{(0, k_2, k_3)}$. For
$k_2=(1-\alpha)/4\alpha,k_3=e^{-d\Delta t}$, $\Sp(F-KC)=\cb{e^{-d\Delta t}, \pp{1\pm e^{-d\Delta t/2}}/2}$,
which fulfills the condition for complete detectability.
Similarly, the pair $[F,NG]$ is completely stabilizable if 
$\exists K\in\R^3: \abs{\lambda}<1, 
\forall\lambda\in\Sp(F-NG\transp{K})$ \cite{anderson_1979}.
Here $G=\sqrt{Q}>0$. Taking $K=\transp{(0,\frac{1}{G(1-\alpha)\Delta t},0)}$,
$\Sp(F-NG\transp{K})=\cb{e^{-d\Delta t},0}$, which proves complete
controllability. Hence $P_\infty$ and $K_\infty$ are well-defined.

From (\ref{eq:kalman_discrete}), we have $x_{k+1}=F x_k + b_k$, therefore 
the estimation error for the asymptotic gain $K_\infty$ is
\begin{equation}\label{eq:kalman_error}
x_{k+1} - \obsx_{k+1|k} = (F - K_\infty C)(x_k - \obsx_{k|k-1}) - K_\infty v_k.
\end{equation}
For $v_k=0$, this converges to zero since $\norm{F-K_\infty C}<1$.

Finally, if $\exists k: x_k-\obsx_{k|k-1}=0$, it follows from
(\ref{eq:kalman_error}) that $x_{k+l}-\obsx_{k+l|k+l-1}=0,\forall l>0$.
This concludes the proof.

\section{Approach for unknown $\mu(\cdot)$}\label{unknown}

The lack of knowledge of the specific growth rate function eliminates the
possibility of classical approaches for state estimation.
To remedy this issue, a Bayesian regularization approach is proposed for the online estimation of the reaction rate as an unknown function of
time $r(t)=\mu(x(t),t)e(t)$.
This is obtained by augmenting the system (\ref{eq:fluo}) with a
probabilistic prior on $r(t)$ in the form of linear stochastic dynamics
\cite{pavlou_2022,denicolao_2003},
and directly solving the augmented linear estimation problem via the hybrid KF.
Consider the following stochastic prior on $r$
\begin{equation}\label{eq:prior}
    \drt(t) = -\theta r(t) + \kappa w(t)
\end{equation}
where $\cb{w(t)}$ is a standard Gaussian noise process.
Together with the assumption that $r(0)$ is Gaussian with mean zero and
variance $\kappa^2/(2\theta)$, this defines a stationary process.
This convenience choice is discussed for a related application in
\cite{pavlou_2022}.
In essence, $\theta$ defines the timescale of admissible fluctuations
(the larger the $\theta$, the faster the fluctuations),
and the standard deviation $\kappa/\sqrt{2\theta}$ defines the admissible
magnitude of $r$. This already enables an educated guess for the choice of the
design parameters $\theta$ and $\kappa$.
Precise tuning of these parameters can be obtained from preliminary data via
the ML procedure discussed in \ref{ML}.

For the ML-optimal values $\obstheta$ and $\obsk$, simultaneous online estimates
of state $x$ and $r$ are obtained from measurements $y_k$ by the KF
for the joint system (\ref{eq:chemostat})--(\ref{eq:prior}),
\begin{equation}\label{eq:augmented}
    \ddt\bmat{x(t)\\r(t)} =
    \bmat{-dI & N\\0 & -\obstheta}\bmat{x(t)\\r(t)} + \bmat{d\Xin\\0}
    + \bmat{0\\\obsk}w(t).
\end{equation}

\begin{figure}[b]
    \includegraphics[width=\columnwidth]{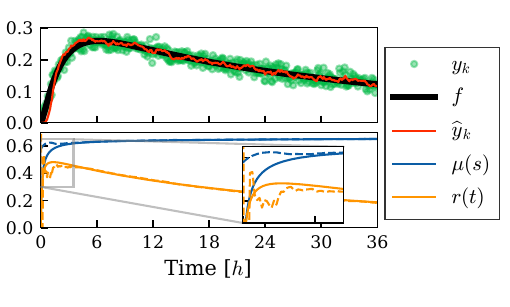}\newline\vspace{-2.4em}
    \caption{Measurements and the reaction rate plotted against time.
    [Top] The fluorescent protein concentration $f(t)$ \textit{(black curve)},
    noisy fluorescence measurements $y_k$ \textit{(green dots)}, and the
    pre-smoothed measurements $\obsy_k$ \textit{(red curve)}.
    [Bottom] The specific growth rate $\mu(s(t))$ and its pre-estimation
    via the conservation law $\mu_k$
    \textit{(blue curves, solid and dotted, respectively)},
    and the total growth rate $r(t)=\mu(s(t)e(t))$ and its BKF estimation
    $\obs{r}_k$ \textit{(orange curves, solid and dotted, respectively)}.}
    \label{fig:measurement}
\end{figure}

Compared with the approach in Sec. \ref{known}, no further noise
component appears in the dynamics of $x$. Performance of this approach
will be demonstrated in simulation in Sec. \ref{results}.\begin{figure*}[t]
    \includegraphics[width=\textwidth]{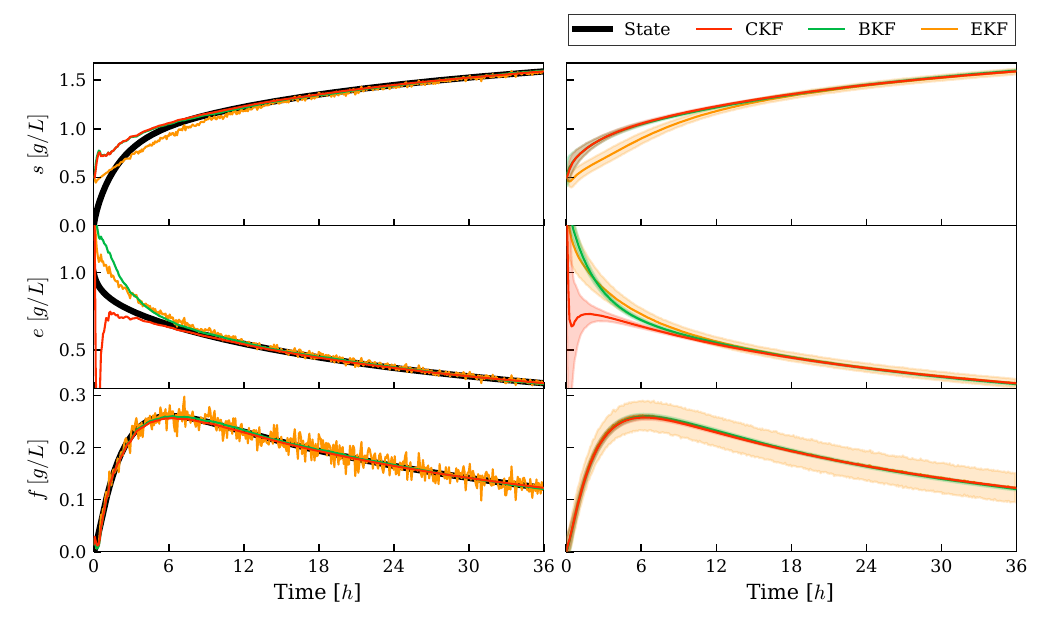}
    \vspace{-1em}\caption{Comparison of KF for state estimation from fluorescence measurements.
    [Left] Reconstructed trajectories using CKF, BKF, and EKF using set of
    generated noisy data compared with the true state trajectory.
    [Right] Mean trajectory of estimated states for each filter over 500 data
    sets and a confidence interval given by twice the standard error around
    the mean trajectory.}
    \label{fig:trajectory}
\end{figure*}

\section{Simulation results}\label{results}

In this section, numerical simulations of the proposed approaches are presented
for system (\ref{eq:fluo}) with Monod growth rate function.
The system is simulated over $36$ hours, with measurements taken at a fixed
interval of $5$ minutes.
Parameters are fixed arbitrarily but realistically based on the
related system in \cite{mauri_2020}.
We use $\gamma=1, k_s=0.2~\gL, \alpha=0.3, \Sin=2~\gL$, the maximal growth
rate $\mumax=\ln2~\h$ corresponds to a doubling time of $1~h$, and 
we fix $d=0.48~\h$.
The initial conditions for the simulation are
$x(0)=\transp{(0,1,0)}~\gL$ whereas $\obsx(0)=\transp{(0.5,1.5,0)}~\gL$.

For the approach presented in Sec. \ref{known}, $\obssig_1(0)=0.5~\gL$
is selected arbitrarily. The ML-optimal smoothing parameter is
$\obs{\beta}\approx 0.013$.
This method being based on the system's conservation
law, we refer to it as the Conservation KF (CKF).
As for the Bayesian regularization approach in Sec. \ref{unknown},
we call it the Bayesian KF (BKF). Its optimal process parameters are
$\obstheta\approx0.03$ and $\obsk\approx0.002$.
The ML problems are solved through the optimizers provided in the Python
library \textit{SciPy} over numerically generated \textit{preliminary data},
with the same initial conditions and the noise covariance $R$.
Different data are later generated for running the filters.

The pre-estimates from the CKF pipeline are shown in Fig. \ref{fig:measurement}.
Estimate $\obsmu_k$ converges to $\mu(s(t))$ at rate determined by $d$ after an
onset resulting from the difference $\sigma_1(0)-\obssig_1(0)$.
The quality of the ML-estimates of $\beta$ (for CKF) and of $\theta$ and
$\kappa$ (for BKF) is witnessed by the estimates $\obsy$, which duly provide
a smoothened version of the data, and the performance of $\obs{r}$,
which converges to $r$ after an onset. 
Reported in \ref{fig:measurement} in the interest of illustration, $\obs{r}$ is obtained
in practice along with the estimation of the other state variables $(e,s,f)$,
which is the focus of the next paragraph.

We then run each filter on a noisy dataset to get estimates $\obsx_k$ of $x(t)$.
For comparison, we also run
the well-known Extended KF (EKF). The $Q$ parameter for this filter
is also fixed by solving the ML problem of subsection \ref{ML} on preliminary
data, with 3 degrees of freedom, yielding
$\obs{Q}\approx\text{diag}(0.04, 0.004, 0.00002)$.
Estimation results on a single dataset are given in Fig. \ref{fig:trajectory} (left).
We also evaluate each filter statistically by Monte Carlo analysis,
\textit{i.e.} by estimating the same $x(t)$ from $N_{MC}=500$ different noisy
data sets with same noise intensity.
A mean estimated trajectory is obtained for each filter as well as a $95\%$ confidence band.
Results from this statistical analysis are in Fig. \ref{fig:trajectory} (right). 
Results show that
all studied filters perform well generally.
The estimated trajectories seem to converge to the true trajectory. However,
both the CKF and BKF converge faster than EKF.

\section{Conclusions}\label{conclusion}

We proposed two online state estimators for a bacterial growth
model in a bioreactor that includes synthesis of a fluorescent reporter protein. The estimators reconstruct the
three-dimensional system's state from discrete-time noisy fluorescence measurements.

We analyzed the system observability showing the possibility to construct a
tunable observer, which we achieved by our CKF.
Using conservation laws, this solution pre-estimates the growth rate profile to
turn state estimation into a linear time-varying problem. This allowed us to
benefit from the KF convergence properties to also prove deterministic
convergence in absence of noise.
Our study focused on the Monod growth function, yet the CKF can be extended to
other growth laws depending on both substrate and biomass concentration,
such as the Contois law \cite{bastin_1990}.
Indeed, for such a growth law, the a-priori unknown biomass concentration
profile can be pre-estimated based on the second conservation relation in
(\ref{eq:cons_fluo}), obtaining $\obse_\sigma$ from the fluorescence
measurements sided by a detector for the conservation variable.
In more generality, the approach is viable for any growth model depending on
state variables that enter suitable conservation laws.

We also proposed what we called BKF as a solution that does not assume knowledge of the
growth rate law, a relevant scenario for bioprocesses. The approach allows estimating the system's state efficiently when classical approaches are not
applicable. By its nature, the BKF can be applied to a more general class of growth
rate functions, provided a Bayesian prior can be fixed to describe the system's
reaction rate behavior over time.


\bibliography{ref.bib}
\bibliographystyle{IEEEtran}

\end{document}

%% file: math.tex
\def\R{\mathbb{R}}

\def\phi{\varphi}

\def\Sin{s_{\text{in}}}
\def\Xin{x_{\text{in}}}

\def\mumax{\mu_{\max}}

\def\dxt{\dot{x}}
\def\dst{\dot{s}}
\def\det{\dot{e}}
\def\dft{\dot{f}}
\def\drt{\dot{r}}

\def\dyt{\dot{y}}

\def\dft{\dot{f}}
\def\dsigt{\dot{\sigma}}

\def\obs#1{\hat{#1}}
\def\err#1{\tilde{#1}}
\def\obsx{\obs{x}}
\def\obsy{\obs{y}}
\def\obss{\obs{s}}
\def\obse{\obs{e}}

\def\obsmu{\obs{\mu}}

\def\obssig{\obs{\sigma}}
\def\obstheta{\obs{\theta}}
\def\obsk{\obs{\kappa}}

\def\erry{\err{y}}
\def\errs{\err{s}}

\def\errmu{\err{\mu}}

\def\dd{\text{d}}

\def\dt{\text{d}t}

\def\ddt{\frac{\dd}{\dt}}

\def\partials#1#2{\frac{\partial#1}{\partial#2}}

\def\pp#1{\left(#1\right)}
\def\bb#1{\left[#1\right]}
\def\cb#1{\left\{#1\right\}}
\def\abs#1{\left\lvert#1\right\rvert}
\def\norm#1{\left\lVert#1\right\rVert}

\def\bmat#1{\begin{bmatrix}#1\end{bmatrix}}
\def\pmat#1{\begin{pmatrix}#1\end{pmatrix}}

\def\transp#1{{#1}^{\top}}

\def\L{L}
\def\gL{g\L^{-1}}
\def\h{h^{-1}}

\def\CC{\mathcal{C}}
\def\OO{\mathcal{O}}
\def\Tk{{\mathcal{T}_k}}

\DeclareMathOperator{\rk}{rank}

\DeclareMathOperator{\Sp}{Sp}
\def\EE#1{\mathbb{E}\cb{#1}}

\DeclareMathOperator*{\argmin}{argmin}
\DeclareMathOperator*{\argmax}{argmax}